\documentclass{article}

\usepackage[T2A]{fontenc}
\usepackage[cp1251]{inputenc}
\usepackage[english,russian]{babel}
\usepackage[tbtags]{amsmath}
\usepackage{amsfonts,amssymb}

\usepackage[left=2.5cm, top=2cm, right=1cm, bottom=20mm]{geometry}

\usepackage[dvipsnames]{xcolor}

\usepackage{mathrsfs}

\usepackage{graphicx}

\numberwithin{equation}{section}

\newtheorem{theorem}{Теорема}

\newtheorem{lemma}{Лемма}

\newtheorem{problem}{Задача}

\newtheorem{proof}{Доказательство}
\newtheorem{remark}{Замечание}

\newcommand{\ff}{\textbf{f}}

\date{	 }

\author{В.\,Н.~Сивкин \\ НИУ ВШЭ \\ sivkin96@yandex.ru \and А.\,А.~Шкаликов \\ МГУ им. Ломоносова \\ ashkaliko@yandex.ru}

\title{Асимптотика спектров задач Дирихле и Дирихле-Неймана для уравнения Штурма-Лиувилля с интегральным возмущением}

\begin{document}

	\maketitle

		\begin{abstract} 
			
			В статье исследуются задачи Дирихле и Дирихле-Неймана для уравнения Штурма-Лиувилля, возмущенного интегральным оператором с ядром свертки. Найдены точные асимптотические формулы  для собственных значений этих задач. Формулы содержат информацию о коэффициентах Фурье потенциала и ядра, а для остаточных членов асимптотики получены оценки, в которых учитывается и скорость убывания с ростом номера собственного значения и скорость убывания при стремлении норм потенциала и ядра к нулю. Формулы являются новыми и в случае оператора Штурма-Лиувилля, когда ядро свертки равно нулю.

			Библиография: 26 названий.

		\end{abstract}
		
		\textbf{Ключевые слова:}
			Оператор Штурма-Лиувилля, интегро-дифференциальный оператор, асимптотические формулы для собственных значений, оператор Харди.
			
			\vskip 3mm
			
		\Large{Asymptotics of the spectra of the Dirichlet and Dirichlet-Neumann problems for the Sturm-Liouville equation with integral perturbation}
		
		\normalsize
		
		\vskip 3mm
		
		\begin{center}
			A.A.Shkalikov, V.N.Sivkin
		\end{center}

		\begin{otherlanguage}{english}
		
			\begin{abstract} 
			
		The article studies the Dirichlet and Dirichlet-Neumann problems for the Sturm-Liouville equation perturbed by an integral operator with a convolution kernel. Sharp asymptotic formulas for the eigenvalues of these problems are found. The formulas contain information about the Fourier coefficients of the potential and the kernel, and estimates are obtained for the remainder terms of the asymptotics, which take into account both the rate of decrease  as the eigenvvalues tend to infinity and the rate of decrease as the norms of the potential and kernel tend to zero. The formulas are also new in the case of the Sturm-Liouville operator, when the convolution kernel is zero.

		Bibliography: 26 items.

		\end{abstract}
		
		\end{otherlanguage}

		\normalsize
		

		\footnotetext{Исследование поддержано Российским научным фондом,
			грант № 25-11-00304.}

		\section{Введение}
		
		В этой статье мы рассматриваем оператор Штурма-Лиувилля, возмущенный интегральным оператором свертки
		\begin{align}\label{eq:LqM}
		L(q, M) = -y'' +q(x)y + \int_0^x M(x-t)y(t)dt.
		\end{align}
		
		Предполагаем, что функция $q \in L_2(0, \pi),$ а функция $M \in L_{2, \pi} := \{\varphi | (\pi-x)\varphi(x) \in L_2(0, \pi)\}, $
		то есть,
		\begin{align}\label{eq:M}
		M(x) = (\pi-x)^{-1} M_0(x), \quad M_0 \in L_2(0, \pi).
		\end{align}
		
		Обе функции $q$ и $M$ предполагаем комплекснозначными. Далее для пространств $L_2(0, \pi)$ и $L_{2, \pi}(0, \pi)$ будем использовать сокращенную запись $L_2$ и $L_{2, \pi}.$ Выбор весового пространства  $L_{2, \pi}$ для функции $M$ не случаен. Он продиктован многочисленными предшествующими работами о рассматриваемом операторе; естественность такого выбора будет ясна также из содержания настоящей статьи.
		
		Далее нам удобно обозначать пару функций $\{q, M_0\}$ в пространстве $L_2 \oplus L_2$ одним символом
		\begin{align}
		\ff = \{q, M_0\}, \quad \|\ff\|_{L_2\oplus L_2} = \|q\|_{L_2} + 2 \|M_0\|_{L_2},
		\end{align}
		(выбор такой нормы будет ясен из дальнейшего).
		
		С оператором  $L(q, M)$ свяжем спектральные задачи, порождаемые краевыми условиями Дирихле и Дирихле-Неймана
		\begin{align}
		& L(q, M)y = \lambda y, \qquad y(0) = y(\pi) = 0, \label{eq:D} \\
		& L(q, M)y = \mu y, \qquad y(0) = y'(\pi) = 0, \label{eq:DN}
		\end{align}
		где $\lambda$ и $\mu$ спектральные параметры.
		
		Большинство работ, связанных с изучением оператора \eqref{eq:LqM}, посвящено исследованию обратных задач для этого оператора. Первой существенной работой была статья В.А.Юрко  \cite{Y1991} (см. также \cite{FY}, \cite{Y2007}). В частности, из результатов \cite{Y1991} следует, что собственные значения $\{\lambda_k\}$ этого оператора подчинены асимптотике
		\begin{align}
		\lambda_k^{1/2} = k + \frac{[q]}{2k} + \frac{\delta_k}{k}, \quad [q] =  \frac{1}{\pi}  \int_0^{\pi} q(t) dt, \quad k=1, 2, ..., \label{eq:AsD}
		\end{align}
		где $\{\delta_k\}_{k=1}^{\infty}$ --- последовательность из пространства $l_2.$ Но основное содержание статьи \cite{Y1991} связано с постановкой обратной задачи, формулируемой следующим образом.
		
		\begin{problem}
			Пусть $q$ --- фиксированная функция. Зная спектр задачи Дирихле \eqref{eq:D}, найти функцию $M.$
		\end{problem}
		
		Работа \cite{Y1991} вызвала большую серию последующих статей  \cite{B2007} -- \cite{B2022}, где эта обратная задача получила развитие. Наиболее полное решение представлено в недавней статье С.А.Бутерина \cite{B2021-2}.
		
		Вторая постановка обратной задачи для оператора \eqref{eq:LqM} была предложена в статье Ю.С.Курышовой \cite{K2007}, которая выполнялась под руководством В.А.Юрко. Задача состоит в следующем:
		\begin{problem}
			Зная функцию $M$ и спектры $\{\lambda_k\}$ и $\{\mu_k\}$ задачи Дирихле \eqref{eq:D} и Дирихле-Неймана \eqref{eq:DN}, найти функцию  $q.$
		\end{problem}
		
		Авторам известны работы Ю.В.Курышовой и С.Т.Шиха \cite{KS}, а также работа В.А.Юрко \cite{Y2017},  где были продолжены исследования этой обратной задачи, но пути полного её решения пока не ясны.
		
		В этой работе мы обратимся к решению прямых задач на собственные значения для \eqref{eq:D} и \eqref{eq:DN}. Во-первых следует отметить, что без труда можно получить асимптотику $\{\mu_k\}$ для задачи \eqref{eq:DN}, аналогичную \eqref{eq:AsD}. А именно,
		\begin{align}
		\mu_k^{1/2} = k - \frac{1}{2} + \frac{[q]}{2k-1} + \frac{\delta_k}{k}, \quad k = 1, 2, ..., \quad \{\delta_k\} \in l_2. \label{eq:AsDN}	
		\end{align}
		Однако асимптотики  \eqref{eq:AsD} и \eqref{eq:AsDN} не дают никакой информации о функциях  $q$ и $M$  (кроме среднего значения $[q]$ функции $q$). Цель этой работы  --- получить более тонкие асимптотики для последовательностей $\{\lambda_k\}$ и $\{\mu_k\},$ в которых информация о функциях $q$ и $M$ уже будет содержаться. Тем самым, нашей целью является не только детальное решение прямых задач для оператора \eqref{eq:LqM}, но и последующее решение некоторых обратных задач, которые в этой статье рассматриваться не будут.
		
		Напомним обозначение $\ff  = \{q, M_0\}$ для функций из пространства $L_2 \oplus L_2.$  Через ${\cal O}(\ff)$ и ${\cal O}(\ff^{\,2})$ будем обозначать величины (функции или числа), зависящие от $\ff$ и, возможно, от индекса $k,$ для которых справедливы оценки
		\begin{align*}
		|{\cal O}(\ff\,)| \leq C\|\ff\,\|, \quad |{\cal O}(\ff^{\,2})| \leq C\|\ff\|^2.
		\end{align*}
		
		В этих оценках предполагаем, что постоянные не зависят ни от индекса $k,$
		ни от $\ff,$ если только $\ff$ изменяется в шаре фиксированного радиуса $r_0.$ В частности, постоянные абсолютные, если $\ff$ лежит в единичном шаре пространства $L_2\oplus L_2.$ \noindent
		В этих оценках $|\cdot|$ означает модуль числа или норму функции в пространстве $C[0, \pi].$
		
		Сформулируем основные результаты работы. В формулировке участвуют оператор Харди $H$ и оператор инволюции $J,$ о которых скажем в следующем параграфе.
		
		\begin{theorem}
			Собственные значения $\{\lambda_k\}$ задачи Дирихле \eqref{eq:D} подчинены асимптотике
			\begin{align}\label{eq:lak}
			\lambda_k^{1/2} = k + \frac{[q]}{2k} - \frac{a_{2k}}{2(2k)} + \frac{m_{2k} - n_k}{2(2k)} + \frac{\delta_k}{k^2}, \quad k = 1, 2, ...,
			\end{align}
			где
			\begin{align*}
			a_k  &= \frac{2}{\pi} \int_0^{\pi} q(t) \cos kt dt, \\
			m_k &= \frac{2}{\pi} \int_0^{\pi} M_0(t) \cos\left( \frac{kt}{2} \right)dt, \\
			n_k &= \frac{2}{\pi} \int_0^{\pi} TM_0(t) \cos kt dt.
			\end{align*}
			Здесь $[q]$ определено в \eqref{eq:AsD}, $T = JH^*J,$ где  $J$ --- оператор инволюции, $H^*$ --- сопряженный к оператору Харди, а последовательность $\{\delta_k\}_{k=1}^{\infty}$ в \eqref{eq:lak} такова, что
			\begin{align}\label{eq:f2}
			\left(\sum_{k=1}^{\infty}|\delta_k|^2\right)^{1/2} = {\cal O}(\ff^{\,2}).
			\end{align}
			
		\end{theorem}

		В этой формулировке теоремы важно, что условие на последовательность $\{\delta_k\}_{k=1}^{\infty},$ которая фигурирует в остаточном члене \eqref{eq:lak} мы получаем не в обычной форме $\{\delta_k\}_{k=1}^{\infty} \in l_2,$ как это проводилось в предшествующих работах,  а в более точной форме \eqref{eq:f2}. Обратим также внимание на то, что правая часть в \eqref{eq:f2} оценивается величиной $\leq C\|\ff\|^2,$ где $C$ --- абсолютная константа, если $\ff\,$ лежит в единичном шаре.
		
		\begin{theorem}
			Спектр задачи Дирихле-Неймана  \eqref{eq:DN} подчинен асимптотике
			\begin{align}\label{eq:muk}
			\mu_k^{1/2} = k - \frac{1}{2} + \frac{[q]}{2k-1} - \frac{a_{2k-1}}{2(2k-1)} + \frac{m_{2k-1}}{2(2k-1)} + \frac{\delta_k}{k^2}, \quad k = 1, 2, ...,
			\end{align}
			где числа $[q],$ $\{a_k\}$ и $\{m_k\}$ определены в Теореме 1, а для последовательности $\{\delta_k\}_1^{\infty}$ выполнена оценка \eqref{eq:f2}.
		\end{theorem}
		
		Для оператора Штурма-Лиувилля без интегрального возмущения эти теоремы были получены в работах
		\cite{SSh2008} -- \cite{SSh2003},  но без  оценки остатка   в форме \eqref{eq:f2}.
		Для доказательства этих теорем нам потребуется предварительный результат о грубой асимптотике собственных значений. Доказательство его проводится с помощью стандартного приема --- теоремы Руше, но все же новые нюансы возникают. Они связаны с тем, что оценки уклонений собственных значений мы проводим
		с константами, не зависящими от функции $\ff$ и индекса $k$, причем не только для  достаточно далеких $k\geq k_1$), но и для близких собственных значений (что необходимо для получения оценок вида \eqref{eq:f2}). Утверждения такого рода мы не встречали в литературе даже для обычного оператора Штурма-Лиувилля (при $M\equiv 0$), поэтому имеет смысл сформулировать его в виде теоремы.
		\begin{theorem}
			Пусть $\lambda_k$ --- собственные значения задачи \eqref{eq:D} и $\rho_k = \lambda_k^{1/2}$ с условием $-\pi/2 < arg\, \rho_k \leq \pi/2.$ Существуют абсолютные постоянные $c_1$ и $c_2,$ такие, что $\forall k \geq 1, \, \forall \ff \in L_2\oplus L_2$
			\begin{align}\label{eq:uni}
			|\rho_k-k| \leq C\|\ff\| k^{-1}, \quad C = \max(c_2 r, c_1),
			\end{align}
			где $r$ --- радиус шара, в котором лежит  $\ff.$ В частности, для всех  $\ff$ из шара $\|\ff\|\leq r_0 = c_1/c_2$ оценка \eqref{eq:uni} выполняется с константой $c_1.$ В качестве константы $c_1$
			можно взять любое число $> 1/\sqrt\pi$.  Но константа $c_2$  зависит от $c_1$  и увеличивается при уменьшении $c_1$.
		\end{theorem}
		
		Аналогичное утверждение имеет место для собственных значений задачи \\ Дирихле-Неймана. Сформулируем 
		этот результат.
		
		ТЕОРЕМА 3'. {\it Пусть $\mu_k$ --- собственные значения задачи \eqref{eq:DN} и $\rho'_k = \mu_k^{1/2}$ с условием $-\pi/2 < arg\, \rho_k \leq \pi/2.$ Существуют абсолютные постоянные $c_1$ и $c_2,$ такие, что $\forall k \geq 1, \, \forall \ff \in L_2\oplus L_2$
			\begin{align}\label{eq:uni'}
			|\rho'_k-k +1/2| \leq C\|\ff\| (k - 1/2)^{-1}, \quad C = \max(c_2 r, c_1),
			\end{align}
			где $r$ --- радиус шара, в котором лежит  $\ff.$   Постоянные  $c_1$ и $c_2$  могут быть взяты такими же, как в Теореме 3.}
		
		Отметим, что числа $(k - 1/2)^{-1}$  в правой части неравенства  \eqref{eq:uni'}  можно заменить на 
		$k^{-1}$, но постоянные в оценках тогда будут другими.
		
		Коротко скажем о дальнейшем плане изложения. В следующем \S\,2 мы высвечиваем полезную роль оператора Харди в рассматриваемой задаче. В частности, доказываем две леммы, которые играют ключевую роль в доказательстве основных теорем. В \S\,3 мы доказываем Теорему 3, а в \S\,4 Теорему 1. В \S\,5 отмечаем изменения, которые надо сделать для доказательства Теоремы 2.
		
		\section{Оператор Харди и его применение для преобразования интегральных выражений и их оценки}
		
		Полезную роль в дальнейшем играют оператор Харди $H$ и оператор инволюции $J,$ определяемые равенствами
		\begin{align*}
		H\varphi(x) = \frac{1}{x} \int_0^x \varphi(t)dt, \quad J\varphi(x) = \varphi(\pi-x).
		\end{align*}
		
		Предполагаем, что эти операторы рассматриваются в пространстве $L_2(0, \pi).$ Очевидно, оператор $J$ самосопряжен и $J^2 = 1.$ Известно \cite{BHS} (см. также \cite{Naz}), что оператор Харди ограничен в $L_2(0, \pi)$ и его сопряженный имеет вид
		\begin{align*}
		H^*\varphi(x) = \int_x^{\pi} \frac{\varphi(t)}{t} dt.
		\end{align*}
		
		Известно также \cite{BHS}, что операторы $1-H$ и $1-H^*$ унитарно эквивалентны операторам левого и правого сдвига в пространстве последовательностей $l_2.$ В частности, $\|H\| = \|H^*\| = 2.$
		
		В дальнейшем важную роль играют следующие две леммы, в которых проясняется роль операторов $H$ и $J.$
		
		\begin{lemma}
			Справедливы следующие равенства:
			\begin{align}
			&
			\begin{aligned}
			&\int_0^x \sin \rho(x-t)\int_0^t M(t-s) \sin (\rho s) ds dt = \\
			&= -\frac{1}{2} \int_0^x (x-z)M(z) \cos \rho(x-z) dz + \frac{1}{2\rho} \int_0^x M(z) \sin \rho(x-z) dz,
			\end{aligned} \label{eq:Main1} \\
			& \begin{aligned}
			&\int_0^{\pi} \sin \rho(\pi-t) \int_0^t M(t-s) \sin \rho s ds dt  = \\
			&= -\frac{1}{2} \int_0^{\pi} M_0(z) \cos \rho(\pi-z) dz + \frac{1}{2} \int_0^{\pi} (TM_0)(z) \cos \rho(\pi-z) dz,
			\end{aligned}\label{eq:Main2}
			\end{align}
			где функции $M$ и $M_0$ связаны равенством \eqref{eq:M}, а оператор $T = JH^*J.$
		\end{lemma}
		
		\begin{proof}
			Преобразуем интеграл в левой части \eqref{eq:Main1}. После замены $t-s = z$ поменяем пределы интегрирования и учтем, что $0\leq z \leq t\leq x.$ Получим
			\begin{align}
			\begin{aligned}
			&2 \int_0^x \sin \rho(x-t) \int_0^t M(t-s) \sin\rho s ds dt = \\
			&= 2 \int_0^x \int_0^t M(z) \sin \rho(x-t) \sin \rho(t-z) dz dt = \\
			&=2 \int_0^x M(z) \int_z^x \sin \rho(x-t) \sin \rho(t-z)dt dz.
			\end{aligned}\label{eq:int1}
			\end{align}
			Внутренний интеграл в правой части этого равенства считаем явно
			\begin{align}
			\begin{aligned}
			&2 \int_z^x \sin \rho(x-t) \sin \rho(t-z) dt dz = \int_z^x (\cos \rho(x+z-2t) - \cos \rho(x-z))dt = \\
			&=-(x-z)\cos \rho(x-z) + \frac{1}{\rho} \sin \rho(x-z).
			\end{aligned}\label{eq:int2}
			\end{align}
			Здесь при переходе ко второму равенству мы воспользовались равенством
			\begin{align}
			\begin{aligned}
			&\int_z^x \cos \rho(x+z-2t) dt = -\frac{1}{2\rho} \sin \rho(x+z-2t)|_{t=z}^{t=x} = \frac{1}{\rho} \sin \rho(x-z).
			\end{aligned}
			\end{align}
			Теперь \eqref{eq:Main1} следует из \eqref{eq:int1} и \eqref{eq:int2}.
			
			Равенство \eqref{eq:Main1} получается после подстановки $x=\pi$ в \eqref{eq:Main1}. Нужно лишь пояснить, что с учетом \eqref{eq:M} и равенства
			\begin{align*}
			\frac{1}{\rho} \sin \rho(\pi-z) = \int_0^{\pi-z} \cos \rho t dt,
			\end{align*}
			второе слагаемое в правой части \eqref{eq:Main1} при $x=\pi$ преобразуется к виду
			
			\begin{align}\label{eq:2.6}
			\begin{aligned}
			&\frac{1}{2\rho} \int_0^{\pi} M(z) \sin \rho(\pi-z) dz = \frac{1}{2} \int_0^{\pi} M_0(z) \left(\frac{1}{\pi-z}\int_0^{\pi-z}\cos \rho t dt\right) dz = \\
			&=\frac{1}{2} \int_0^{\pi} M_0(z) (JHJ)(\cos \rho(\pi-t))(z)dz = \frac{1}{2} \int_0^{\pi} (TM_0)(z) \cos \rho(\pi-z)dz.
			\end{aligned}
			\end{align}

			Лемма доказана.
			
		\end{proof}

		\begin{lemma}
			Пусть функция  $\varphi \in C[0, \pi],$ $M_0 \in L_2(0, \pi),$ а функция $M$ определена равенством \eqref{eq:M}. Тогда функция
			\begin{align*}
			F(t) = \int_0^t M(t-\xi) \varphi(\xi) d\xi \in L_2(0, \pi), 
			\end{align*}
			причем
			\begin{align}
			\|F\|_{L_2} \leq 2\|\varphi\|_C \|M_0\|_{L_2}. \label{eq:base}
			\end{align}
			Аналогичное утверждение справедливо и для функции $F_1(t) = \int_t^{\pi} M(\xi -t)\varphi(\xi) d\xi.$
		\end{lemma}
		
		\begin{proof}
			Имеем
			\begin{align*}
			|F(t)|  &= \left|\int_0^t \frac{M_0(s)}{\pi-s} \varphi(t-s) ds \right| = \left|\int_{\pi -t}^{\pi} \frac{M_0(\pi-\xi)}{\xi} \varphi(t+\xi-\pi) d\xi\right| \\
			& \leq \|\varphi\|_C \left|(JH^*J)(|M_0(\xi)|)(t)\right|. \\
			\end{align*}
			Из оценки $\|JH^*J\|\leq 2$ получаем \eqref{eq:base}.

			Доказательство для функции $F_1$ сохраняется.
		\end{proof}

		\section{Доказательство теоремы 3.}
		
		Обозначим через $S(x, \rho)$ решение уравнения \eqref{eq:D} при $\lambda = \rho^2,$ подчиненное условию $S(0, \rho) = 0,$ $S'(0, \rho) = 1.$ Ясно, что все собственные значения задачи \eqref{eq:D} определяются нулями $\rho_k = \lambda_k^{1/2}$ целой функции $S(\pi, \rho).$ Эта функция четная, поэтому далее рассматриваем только нули, подчиненные условию $-\pi/2<arg\, \rho_k\leq \pi/2,$ которые полностью определяют собственные значения $\lambda_k.$ Нумерацию ведем с учетом алгебраической кратности в порядке возрастания модулей.
		
		Функция $S(x, \rho)$ удовлетворяет интегральному уравнению
		\begin{align}
		S(x, \rho) = \frac{\sin \rho x}{\rho} + \int_0^x \frac{\sin \rho(x-t)}{\rho} \left(q(t)S(t, \rho)+\int_0^t M(t-s)S(s, \rho)ds\right)dt,
		\end{align}
		\noindent
		эквивалентному уравнению \eqref{eq:D} при выбранных начальных условиях. Решая это уравнение методом последовательных приближений, представим его решение в виде ряда
		\begin{align*}
		S(x, \rho) = \sum_{n=0}^{\infty} S_n(x, \rho),
		\end{align*}
		где
		\begin{align}
		\begin{aligned}
		&S_0(x, \rho) = \frac{\sin \rho x}{\rho}, \quad S_0(x, 0) = x, \\
		&S_n(x, \rho) = \int_0^x \frac{\sin \rho(x-t)}{\rho}\left(q(t)S_{n-1}(t, \rho) + \int_0^t M(t-s)S_{n-1}(s, \rho)ds\right)dt.
		\end{aligned} \label{eq:Sind}
		\end{align}
		
		\begin{lemma}
			Справедливы оценки
			\begin{align}\label{eq:Sn}
			|S_n(x, \rho)| \leq \frac{\|\ff\,\|^n e^{\nu x} x^{n/2}}{\sqrt{n!}|\rho|^{n+1}}, \quad 0\leq x \leq \pi, \quad \nu = |Im \rho|, \quad n=0, 1, ... \, .
			\end{align}
		\end{lemma}
		
		\begin{proof}
			При  $n=0$ это утверждение верно. В предположении, что оно верно при $n-1,$ получаем (с учетом $|\sin\rho(x-t)|\leq e^{\nu(x-t)}$)
			\begin{align}\label{eq:Sn0}
			\begin{aligned}
			|S_n(x, \rho)| \leq &  \frac{\|\ff\,\|^{n-1}}{|\rho|\sqrt{(n-1)!}|\rho|^{n}}  \int_0^x e^{\nu(x-t)}\left(|q(t)|t^{(n-1)/2}e^{\nu t} + \right. \\
			+&\left. \int_0^t |M(t-s)|e^{\nu s} s^{(n-1)/2}ds\right)dt.
			\end{aligned}
			\end{align}
			Здесь первое слагаемое (без учета внеинтегрального множителя) в силу неравенства Коши-Буняковского допускает оценку
			\begin{align}\label{eq:Sn1}
			\leq e^{\nu x} \|q\| \left(\int_0^x t^{n-1}dt\right)^{1/2} = \frac{e^{\nu x} \|q\| x^{n/2}}{\sqrt{n}}.
			\end{align}
			Второе слагаемое оценивается так:
			\begin{align}\label{eq:Sn2}
			\begin{aligned}
			&\leq \int_0^x e^{\nu(x-t)}e^{\nu t} t^{(n-1)/2} \int_0^t \frac{|M_0(t-s)|}{\pi-t+s} ds dt  =\\
			&=e^{\nu x} \int_0^x t^{(n-1)/2} \int_{\pi-t}^{\pi} \frac{|M_0(\pi-s)|}{s} ds dt = \\
			&=e^{\nu x } \int_0^x t^{(n-1)/2}(JH^*J)(|M_0(s)|)(t) dt \leq e^{\nu x} \frac{x^{n/2}}{\sqrt{n}} 2\|M_0\|.
			\end{aligned}
			\end{align}
			
			Здесь при переходе к последнему неравенству мы воспользовались неравенством Коши-Буняковского и равенством $\|JH^*J\| = 2.$

			Вспоминая о множителе перед интегралом \eqref{eq:Sn0} и учитывая равенство $\|\ff\| = \|q\|+2\|M_0\|,$ после суммирования \eqref{eq:Sn1} и \eqref{eq:Sn2} приходим к оценке \eqref{eq:Sn}. Лемма доказана.

		\end{proof}
		
		Функцию $S(\pi, \rho)$ представим в виде
		\begin{align*}
		S(\pi, \rho) = \frac{\sin \rho \pi}{\rho} + \frac{1}{\rho} \Phi(\rho).
		\end{align*}
		
		Из этого представления следует, что функция $\rho^{-1}\Phi(\rho)$ целая. Из Леммы 3 следует, что ряд, представляющий эту функцию, сходится. Кроме того, при $|\rho|>\sqrt{\frac{\pi}{2}}\|\ff\|$ справедлива оценка
		\begin{align}
		\begin{aligned}\label{eq:Phi}
		&|\Phi(\rho)| \leq |\rho|\sum_{k=1}^{\infty} |S_k(\pi, \rho)| \leq e^{\nu \pi}\left(t+\frac{t^2}{\sqrt{2!}}+\frac{t^3}{\sqrt{3!}}+...\right) < \\
		&<e^{\nu \pi}t\left(1+\frac{t}{\sqrt{2}}+\left(\frac{t}{\sqrt{2}}\right)^2+...\right) = \frac{e^{\nu \pi} \sqrt{2}t}{\sqrt{2}-t}, \quad \nu = |Im \rho|, \quad t = \frac{\sqrt{\pi}\|\ff\|}{|\rho|}.
		\end{aligned}
		\end{align}
		
		Функция $t(\sqrt{2}-t)^{-1}$ возрастает при $t\in(0, \sqrt{2}),$ поэтому максимум $\Phi(\rho)$ на окружности $|\rho-k| = r \leq 1/2$  (при фиксированном $k\in \mathbb{N}$) достигается при $\rho = k-r.$ Следовательно, на таких окружностях верна оценка
		\begin{align*}
		|\Phi(\rho)| < \frac{e^{\pi r} \sqrt{2\pi} \|\ff\|}{(k-r)\left(\sqrt{2}- \frac{\sqrt{\pi}\|\ff\|}{k-r}\right)} = \frac{e^{\pi r} \sqrt{\pi} \|\ff\|}{k-r- \sqrt{\frac{\pi}{2}} \|\ff\|}, \quad 0<r\leq \frac{1}{2}.
		\end{align*}

		Функция  $\sin \rho \pi$ на таких окружностях оценивается снизу:
		\begin{align*}
		&\min_{|\rho-k|=r} |\sin \rho \pi|  = \min_{|\rho-k|=r} |\sin (\rho-k) \pi|  \geq \\
		&\geq \pi r \left(1-\frac{\pi^2 r^2}{6}\left(1+\frac{\pi^2 r^2}{4 \cdot 5} + \left(\frac{\pi^2 r^2}{4\cdot 5}\right)^2 + ...\right)\right) = \pi r \left(1 - \frac{10\pi^2 r^2}{3(20-\pi^2r^2)}\right).
		\end{align*}
		
		Положим $\alpha = \frac{1}{\sqrt{\pi}}(1+\varepsilon),$ а число $\varepsilon\in (0, 1]$ выберем позже. Обозначим через $\gamma_k$ окружности
		\begin{align*}
		|\rho-k| = \alpha \|\ff\,\|k^{-1} =:r_k.
		\end{align*}
		
		На окружности $\gamma_k$ справедливы оценки (учитываем, что $\alpha^{-1}<\sqrt{\pi}$):
		\begin{align*}
		|\Phi(\rho)| < \frac{e^{\pi r_k} \sqrt{\pi} r_k k \alpha^{-1}}{k - r_k - \sqrt{\frac{\pi}{2}} r_k k \alpha^{-1}} = \frac{e^{\pi r_k} \sqrt{\pi} \alpha^{-1} r_k}{1-(\frac{1}{k}+\frac{\pi}{\sqrt{2}})r_k}.
		\end{align*}

		В силу теоремы Руше внутри окружности  $\gamma_k$ находится в точности один ноль $\rho_k$ функции $S(\pi, \rho),$ если при $r = r_k$ выполняется неравенство
		\begin{align*}
		\pi r \left(1 - \frac{10 \pi^2 r^2}{3(20-\pi^2 r^2)}\right) > \frac{e^{\pi r} \sqrt{\pi} r \alpha^{-1}}{1-(\frac{1}{k}+\frac{\pi}{\sqrt{2}})r }.
		\end{align*}
		
		Это неравенство эквивалентно следующему
		\begin{align}\label{eq:fin}
		(1+\varepsilon)\left(1- \frac{10 \pi^2 r^2}{3(20-\pi^2 r^2)}\right) > \frac{e^{\pi r}}{1-(\frac{1}{k}+\frac{\pi}{\sqrt{2}})r}, \quad k \geq 1.
		\end{align}

		Функция аргумента $r$ в левой части этого неравенства возрастает при $r\to 0,$ а в правой части убывает. При $r=0$ это неравенство выполнено при любом $\varepsilon>0.$ Поэтому при любом $\varepsilon>0$ найдется число $r_0 = r_0(\varepsilon)>0,$ такое, что \eqref{eq:fin} выполняется при всех $r = r_k \leq r_0.$ В частности, легко проверить, что при  $\varepsilon = 1/3$ неравенство выполняется при всех $r = r_k\leq \varepsilon/8  = 1/24,$ $k\geq 1.$ Тем самым, если
		\begin{align}\label{eq:c0}
		\alpha \|\ff\| k^{-1} \leq 1/24, \quad \text{или} \quad k\geq \frac{24}{\sqrt{\pi}}\left(1+\frac{1}{3}\right)\|\ff\| = c_0 \|\ff\|, \quad c_0 = \frac{32}{\sqrt{\pi}},
		\end{align}
		то внутри окружности $\gamma_k$ лежит в точности один нуль $\rho_k,$ причем
		\begin{align}\label{eq:c1}
		|\rho_k - k|< \frac{1}{\sqrt{\pi}}\left(1+\frac{1}{3}\right) \|\ff\|k^{-1} = c_1 \|\ff\|k^{-1}, \quad c_1 = \frac{4}{3\sqrt{\pi}}.
		\end{align}
		
		Пусть $k_1$ --- наименьшее число, при котором выполняется неравенство \eqref{eq:c0}. Покажем, что число остальных нулей $\rho_k,$ подчиненных условию $-\pi/2 < arg \, \rho_k \leq \pi/2,$ равно $k_1-1,$ то есть нумерация этих нулей начинается с индекса $k=1.$
		
		Рассмотрим квадрат $\Pi_k$ с центром в нуле, вертикальные стороны которого проходят через точки $\pm(k+1/2).$ На вертикальных и горизонтальных сторонах прямоугольника  $\Pi_k$ имеем оценки
		\begin{align*}
		&|\sin \rho \pi| = |\sin (\pm(k+1/2)+i\nu)\pi| = |\sin(k+1/2)\pi \cos i\nu \pi| = \frac{e^{\nu \pi}+e^{-\nu \pi}}{2} > \frac{e^{\nu \pi}}{3}, \\
		&|\sin \rho \pi| \geq \frac{e^{\nu \pi} - e^{-\nu \pi}}{2} > \frac{ e^{\nu \pi} }{3}, \quad \nu = |Im \rho|
		\end{align*}
		(второе неравенство на горизонтальных сторонах выполняется при $\nu  \geq 1/(2\pi),$ тем более, при $\nu \geq k+1/2$). Для функции $\Phi$ на сторонах этого квадрата имеем оценку
		\begin{align}
		|\Phi(\rho)| < \frac{e^{\nu \pi} \sqrt{2}t}{\sqrt{2} - t} \leq \frac{e^{\nu \pi} \sqrt{2\pi} \|\ff\|}{\sqrt{2} (k+1/2) - \sqrt{\pi} \|\ff\|  }, \quad t = \frac{\sqrt{\pi} \|\ff\|}{|\rho|}.
		\end{align}
		
		Следовательно, при выполнении неравенства
		\begin{align}\label{eq:k0}
		\frac{1}{3} \geq \frac{ \sqrt{2\pi} \|\ff\|}{\sqrt{2} (k+1/2) - \sqrt{\pi} \|\ff\|  }, \quad \text{или} \quad k \geq c_3 \|\ff\| -\frac{1}{2}, \quad c_3 = (3+ \frac{1}{\sqrt{2}})\sqrt{\pi},
		\end{align}
		\noindent
		функция $S(\pi, \rho)$ имеет внутри квадрата $\Pi_k$ столько же нулей, сколько их имеет функция $\rho^{-1} \sin \rho \pi,$ то есть, ровно $k$ нулей, не считая нули с противоположными знаками.
		
		Пусть $k_0\geq 1$ --- наименьшее число, при котором выполнена оценка \eqref{eq:k0}. Легко видеть, что $k_0\leq k_1,$ причем $k_0 = k_1,$ только если $k_1=1.$ Действительно, если $k_1>1,$ то из определения чисел $k_0$ и $k_1$ следует $c_0\|\ff\|>1,$ а тогда
		\begin{align*}
		k_1-k_0 \geq c_0\|\ff\| - (c_3\|\ff\|-1/2+1) \geq (c_0-c_3 -c_0/2)\|\ff\| > 0.
		\end{align*}
		
		Таким образом, если $k_1 = 1$ (или $\|\ff\|\leq \sqrt{\pi}/32$), то все нули $\rho_k$ подчинены неравенствам \eqref{eq:c1} при $k\geq 1.$ Если же $k_1\geq 2,$ то при $k\geq k_1$ выполнены неравенства \eqref{eq:c1}, а остальные $k_1-1$ нулей лежат в правой полуплоскости внутри $\Pi_{k_1-1},$ так как $k_0 \leq k_1-1.$ Диагональ 1/4 квадрата $\Pi_{k_1-1}$ равна $\sqrt{2}(k_1-1+1/2),$ а из определения числа $k_1$ следует $1\leq k_1-1\leq c_0\|\ff\|.$ Поэтому для нулей $\rho_k \in \Pi_{k_1-1}$ имеем оценки
		\begin{align}\label{eq:c3}
		\begin{aligned}
		|\rho_k-k| < \sqrt{2}(k_1-1+1/2) &\leq (\sqrt{2}(1+1/2)c_0\|\ff\|k)k^{-1} \leq \\
		&\leq \frac{3}{\sqrt{2}} c_0(k_1-1)\|\ff\|k^{-1} \leq c_2 \|\ff\|^2 k^{-1}, \quad c_2 = \frac{3}{\sqrt{2}}c_0^2.
		\end{aligned}
		\end{align}
		
		\noindent Тем самым,
		\begin{align}\label{eq:fin1}
		|\rho_k-k| < C\|\ff\|k^{-1}, \quad  C= \max(c_2 r, c_1), \quad \forall k \geq 1,
		\end{align}
		где $r$ --- радиус шара, в котором находится функция $\ff.$ В частности, если $\ff$ лежит в единичном шаре, то оценка \eqref{eq:fin1} выполнена с абсолютной константой $C.$
		
		Теорема доказана.

		\begin{remark}
			Мы отмечали, что выбор числа $\varepsilon>0,$ которое отвечает за константу $c_1$ в оценке дальних нулей $\rho_k,$ находится в нашей власти. Но уменьшая $\varepsilon$ (или $c_1$), мы увеличиваем константу $c_0$ (то есть, ухудшаем оценку близких нулей) и наоборот. Это полезно иметь в виду при оценке далеких или близких собственных значений.
		\end{remark}
		
		\section{Доказательство Теоремы 1.}
		
		Далее нам удобно будет использовать обозначения
		\begin{equation}\label{eq:abre}
		\{k^{-1}\}:= {\cal O}(\ff)k^{-1}, \quad \{k^{-2}\} := {\cal O}(\ff^{\,2})k^{-2}.
		\end{equation}
		Тогда из  \eqref{eq:fin1} получаем
		\begin{equation*}
		\rho_k = k +\{k^{-1}\}.
		\end{equation*}
		При этом, будем иметь в виду, что  $\{k^{-1}\}$ и  $\{k^{-2}\}$ обозначают числа, не зависящие от переменной $x\in[0, \pi].$ Обозначение $[1]$ будем использовать для функций, зависящих от переменной $x\in [0, \pi]$ и индекса $k\in \mathbb{N},$ модуль которых оценивается либо абсолютной константой $C,$ либо константой $Cr,$ где $r$ --- радиус шара, в котором лежит функция $\ff.$ В случае, если символ $[1]$ стоит перед знаком интеграла, предполагаем, что это число $\leq C$ или $\leq Cr.$ В этих обозначениях имеем представления:
		\begin{align}
		&\cos \rho_k t = \cos kt \cos t\{k^{-1}\} - \sin kt \sin t\{k^{-1}\} = \cos kt +\{k^{-1}\}t \sin kt + [1]\{k^{-2}\}, \label{eq:cos}\\
		&\sin \rho_k t = \sin kt \cos t\{k^{-1}\} + \cos kt \sin t\{k^{-1}\} = \sin kt + \{k^{-1}\} t \cos kt + [1]\{k^{-2}\}, \label{eq:sin}\\
		&\frac{1}{\rho_k} = \frac{1}{k+\{k^{-1}\}} = \frac{1}{k(1+k^{-1}\{k^{-1}\})} = \frac{1}{k} + \frac{1}{k^2} \{k^{-1}\}. \label{eq:rho}
		\end{align}
		
		Отметим, что представления \eqref{eq:cos} и \eqref{eq:sin} справедливы при $k\geq 1,$ а представление \eqref{eq:rho} только для далеких собственных значений при $k\geq k_1,$ где $k_1$ --- наименьшее число, при котором выполняется оценки \eqref{eq:c1}. Числа $\rho_k$ могут обращаться в нуль или быть сколь угодно близкими к нулю, даже если функция  $\ff$ лежит в шаре фиксированного  (но достаточно большого) радиуса. Тем не менее, утверждение Теоремы 1 сохраняется для всех $\ff \in L_2\oplus L_2$ с единственной оговоркой, что оценка ${\cal O}(\ff^{\,2})\leq C\|\ff\|^2$ выполняется с константой $C,$ зависящей от $r = \|\ff\|.$ При $r \leq \sqrt{\pi}/32$ (см. \S 3) константа в этой оценке от $r$  не зависит, число $k_1$ равно $1,$ неравенство $|\rho_k-k|<1/2$ выполняется при всех $k\geq 1,$ поэтому представление \eqref{eq:rho} выполняется при всех $k\geq 1.$  Далее для простоты мы будем считать выполненным условие $\|\ff\| \leq \sqrt{\pi}/32,$ а в конце доказательства теоремы скажем об изменениях, которые нужно сделать в общем случае.
		
		В доказательстве будем использовать следующие два простых предложения.
		\begin{lemma}
			Пусть функция $g(t, s)$ такова, что $|g(t, s)|\leq G(s) \in L_2(0, \pi).$ Положим
			\begin{align}\label{eq:alp}
			\alpha(\rho_k, t) = \int_0^t g(t, s) \cos \rho_k(t-s) ds, \quad t\in [0, \pi],
			\end{align}
			предполагая, что для последовательности $\{\rho_k\}$ справедлива оценка \eqref{eq:rho}. Тогда
			\begin{align}\label{eq:G}
			\sum_{k=1}^{\infty} |\alpha(\rho_k, t)|^2 \leq [1]\|G\|,
			\end{align}
			где постоянная $[1]$ от $t$ не зависит. Это утверждение сохраняет силу, если вмести $\cos \rho_k(t-s)$ в \eqref{eq:alp} фигурирует $\cos \rho_k(t-2s).$
		\end{lemma}
		
		\begin{proof}
			Воспользуемся равенствами \eqref{eq:cos} и \eqref{eq:rho}. Тогда
			\begin{align}
			\alpha(\rho_k, t) = \int_0^t g(t, s) \cos k(t-s) ds  + \{k^{-1}\} \int_0^t [1] G(s) ds.
			\end{align}
			Очевидно, $l_2$-норма последовательности второго слагаемого $\leq [1] \|G\|.$ Первое слагаемое запишем в виде
			\begin{align}
			\cos kt \int_0^{\pi} g(t, s) \chi_t(s) \cos ks ds + \sin kt \int_0^{\pi} g(t, s) \chi_t(s) \sin ks ds = \alpha_1(k, t) +\alpha_2(k, t),
			\end{align}
			где $\chi_t(s)$ --- характеристическая функция отрезка $[0, t].$ Системы $\{\frac{2}{\pi}\cos ks\}_1^{\infty}$ и \\ $\{\frac{2}{\pi}\sin ks\}_1^{\infty}$ ортонормированы в $L_2(0, \pi),$ поэтому из неравенства Бесселя получаем при всех $t\in [0, \pi]$
			\begin{align}
			\sum_{k=1}^{\infty} |\alpha_1(k, t)|^2 + \sum_{k=1}^{\infty} |\alpha_2(k, t)|^2 \leq 2 \|G\|,
			\end{align}
			так как $|g(t, s)\chi_t(s)| \leq G(s).$ Этим завершается доказательство. Оно не меняется при замене $\cos \rho_k(t-s)$ на $\cos \rho_k(t-2s).$ Лемма доказана.
			
		\end{proof}
		
		\begin{lemma}
			Пусть функция $p(t, \rho)$ такова, что $|p
			(t, \rho)| \leq P(t) \in L_2(0, \pi).$ Если
			$$J(\rho_k) = \int_0^{\pi} p(t, \rho_k)\alpha(\rho_k, t) dt,$$
			и для последовательности $\{\alpha(\rho_k, t)\}_1^{\infty}$ выполняется оценка \eqref{eq:G}, то $\{J(\rho_k)\}_1^{\infty} \in l_2$ с оценкой $l_2$-нормы $\leq [1]\|P\|\|G\|.$
		\end{lemma}
		\begin{proof}
			Из неравенства Коши-Буняковского получаем
			\begin{align*}
			&\sum_{k=1}^{\infty} |J(\rho_k)|^2 \leq \sum_{k=1}^{\infty} \int_0^{\pi} |P(t)|^2 dt \int_0^{\pi} |\alpha(\rho_k, t)|^2 dt \leq \\
			&\leq  \|P\| \sum_{k=1}^{\infty} \int_0^{\pi} |\alpha(\rho_k, t)|^2 dt = \|P\| \int_0^{\pi} \sum_{k=1}^{\infty} |\alpha(\rho_k, t)|^2 dt \leq [1]\|P\|\|G\|.	
			\end{align*}
		\end{proof}

		Приступим к доказательству Теоремы 1. Разобьем его на 5 шагов. На первом шаге мы получим асимптотическое  представление последовательности $\{S_1(\pi, \rho_k)\}_1^{\infty}.$ На втором и третьем шагах проведем  оценку норм последовательностей $\{S_2(\pi, \rho_k)\}_1^{\infty}$ и $\{\sum_{n\geq 3} S_n(\pi, \rho_k)\}_1^{\infty}.$ Используя полученные оценки, на четвертом шаге мы получим утверждение Теоремы 1, но при ранее оговоренном условии, когда $\ff$ лежит в круге достаточно малого радиуса $r_0$ (например, $r_0 = \sqrt{\pi}/32$). Наконец, на пятом шаге мы скажем об изменениях, которые нужно провести в общем случае.
		
		\textbf{Шаг 1.} Вместо $S_1(x, \rho)$ нам будет удобнее работать с функцией
		\begin{align}\label{eq:S1}
		\rho^2 S_1(t, \rho) = \int_0^t \sin \rho(t-\xi) \left(q(\xi) \sin \rho \xi + \int_0^{\xi} M(\xi-s) \sin \rho s ds\right)d\xi.
		\end{align}
		
		Преобразуем это выражение, воспользовавшись равенством
		\begin{align}
		&\sin \rho(t-\xi) \sin \rho \xi = \frac{1}{2}(-\cos \rho t + \cos \rho t \cos 2\rho \xi + \sin \rho t \sin 2\rho \xi).
		\end{align}
		
		Согласно этому равенству, а также равенству \eqref{eq:Main1}, выражение \eqref{eq:S1}  перепишем в виде
		\begin{align}
		\begin{aligned}\label{eq:S1t}
		\rho^2 S_1(t, \rho) = & -\frac{1}{2} \cos \rho t \int_0^t q(\xi) d\xi + \frac{1}{2} \cos \rho t \int_0^t q(\xi) \cos 2\rho \xi d\xi + \\
		&+ \frac{1}{2}\sin \rho t \int_0^t q(\xi) \sin 2\rho \xi d\xi - \frac{1}{2} \int_0^t (t-\xi)M(\xi) \cos \rho(t-\xi) d\xi \\
		& + \frac{1}{2\rho} \int_0^t M(\xi) \sin \rho(t-\xi) d\xi.
		\end{aligned}
		\end{align}

		При $t=\pi$ эту функцию представим в виде пяти слагаемых, участвующих в правой части \eqref{eq:S1t}
		\begin{align}\label{eq:I}
		I(\rho) := \rho^2 S_1(\pi, \rho) = \sum_{j=1}^5 I_j(\rho).
		\end{align}
		С учетом равенств \eqref{eq:cos} и \eqref{eq:sin} первые три слагаемых при $\rho = \rho_k$ предстанут в виде
		\begin{align*}
		I_1 &= -\frac{1}{2} \cos k \pi \int_0^{\pi} q(\xi)d\xi +\{k^{-2}\} = (-1)^{k+1}\frac{\pi}{2}[q] + \{k^{-2}\}, \\
		I_2 &= \frac{1}{2} \cos k \pi \int_0^{\pi} q(\xi) \cos 2k\xi d\xi + \{k^{-1}\} \int_0^{\pi} \xi q(\xi) \sin k\xi d\xi + \{k^{-2}\} = \\
		&= \frac{1}{2}(-1)^k \frac{\pi}{2}a_{2k} +\{k^{-1}\}\alpha_k, \\
		I_3 &= \frac{1}{2}\{k^{-1}\} \int_0^{\pi} \xi q(\xi) \cos 2k\xi d\xi + \{k^{-2}\} = \{k^{-1}\} \alpha_k,
		\end{align*}
		где последовательности $\{\alpha_k\}_1^{\infty}$ (они разные в представлениях $I_2$ и $I_3$), таковы, что
		\begin{align}\label{eq:alphak}
		\|\{\alpha_k\} \|_{l_2} = {\cal O}(\ff).
		\end{align}
		
		Четвертое и пятое слагаемые в \eqref{eq:I}, согласно равенству \eqref{eq:2.6}, запишем в виде
		\begin{align*}
		I_4 &=-\frac{1}{2} \int_0^{\pi} M_0(\xi) \cos \rho_k(\pi-\xi) d\xi  = -\frac{1}{2} \int_0^{\pi} M_0(\xi) \cos k(\pi-\xi) d\xi +\\
		&+\{k^{-1}\} \int_0^{\pi} JM_0(\xi) \xi \sin k\xi d\xi +\{k^{-2}\} = \frac{1}{2}(-1)^{k+1} \frac{\pi}{2}m_{2k} +\{k^{-1}\}\alpha_k, \\
		I_5 &= \frac{1}{2} \int_0^{\pi} (TM_0)(z)\cos k(\pi-z) dz +\{k^{-1}\} \int_0^{\pi} (JTM_0)(\xi) \xi \sin k\xi d\xi +\{k^{-2}\} = \\
		&=\frac{1}{2}(-1)^{k} \frac{\pi}{2} n_k +\{k^{-1}\}\alpha_k,
		\end{align*}
		где, как и ранее, последовательности $\{\alpha_k\}_1^{\infty}$ подчинены оценке \eqref{eq:alphak}.
		
		Складывая все 5 слагаемых получаем
		\begin{align}\label{eq:asy}
		I(\rho_k) = (-1)^{k+1}\pi\left(\frac{[q]}{2} - \frac{a_{2k}}{4} + \frac{m_{2k}}{4} -\frac{n_{k}}{4}\right) +\{k^{-1}\} \alpha_k,
		\end{align}
		где последовательность $\{\alpha_k\}_1^{\infty}$ такова, что её $l_2$-норма подчинена оценке \eqref{eq:alphak}.

		\textbf{Шаг 2.} Рассмотрим функцию
		\begin{align}\label{eq:J}
		J(\rho) = \rho^3 S_2(\pi, \rho).
		\end{align}
		Согласно равенствам  \eqref{eq:Sind} и \eqref{eq:S1} эта функция представима в виде четырех слагаемых:
		\begin{align*}
		&J_1(\rho)  = \int_0^{\pi} \sin \rho(\pi-t) q(t) \int_0^t \sin \rho(t-s) q(s) \sin \rho s ds dt, \\
		&J_2(\rho) = \int_0^{\pi} \sin \rho(\pi-t) q(t) \int_0^t \sin \rho(t-s) \int_0^s M(s-\xi) \sin \rho \xi d\xi ds dt,  \\
		&J_3(\rho) = \int_0^{\pi} \sin \rho(\pi-t) \int_0^t M(t-s)  \int_0^s \sin \rho(t-\xi) q(\xi) \sin \rho \xi d\xi ds dt, \\
		&J_4(\rho) = \int_0^{\pi} \sin \rho(\pi-t) \int_0^t M(t-s) \int_0^s \sin \rho(s-\xi) \int_0^{\xi} M(\xi-\tau) \sin \rho \tau d\tau d\xi ds dt.
		\end{align*}

		Покажем, что каждая из последовательностей $\{J_j(\rho_k)\}_1^{\infty}$ принадлежит $l_2$ с оценкой $l_2$-нормы
		\begin{align}\label{eq:Jj}
		\|\{J_j(\rho_k)\}\|_{l_2} \leq [1] \|\ff\|^2, \quad j=1, 2, 3, 4.
		\end{align}

		\textbf{Шаг 2.1.} Внутренний интеграл в выражении для $J_1(\rho)$ представим в виде
		\begin{align*}
		\alpha(\rho, t):&= \int_0^t q(s) \sin \rho(t-s) \sin \rho s ds = - \frac{1}{2} \cos \rho t \int_0^t q(s) ds + \\
		&+ \frac{1}{2} \int_0^t q(s) \cos \rho(t-2s) ds = \alpha_1(\rho, t) + \alpha_2(\rho, t).
		\end{align*}
		Тогда
		\begin{align}
		J_1(\rho) = \int_0^{\pi} q(t) \sin \rho(\pi-t) \alpha_1(\rho, t) + \int_0^{\pi} q(t) \sin \rho (\pi-t) \alpha_2(\rho, t) dt = J_{1, 1}(\rho) + J_{1, 2}(\rho).
		\end{align}
		
		Из Леммы 4 следует, что $l_2$-норма последовательности $\{\alpha_2(\rho_k, t)\}_1^{\infty}$ допускает оценку $\leq [1] \|q\|.$ Тогда из Леммы 5 получаем
		\begin{align}\label{eq:J12}
		\|\{J_{1, 2}(\rho_k)\}\|_{l_2} \leq [1] \|q\|^2 \leq [1] \|\ff\|^2.
		\end{align}
		
		Положим $Q(t) = q(t) \int_0^t q(s) ds.$ Тогда
		\begin{align*}
		&2J_{1, 1}(\rho) = - \int_0^{\pi} Q(t) \sin \rho(\pi-t) \cos \rho t dt = \\ &= - \sin \rho \pi \left(\int_0^{\pi} Q(t)ds + \int_0^{\pi} Q(t) \cos 2\rho t dt \right) +\cos \rho \pi \int_0^{\pi} Q(t) \sin 2\rho t dt.
		\end{align*}
		Пользуясь представлениями  \eqref{eq:cos}, \eqref{eq:sin} и \eqref{eq:rho} из этого равенства получаем
		\begin{align}
		\|\{J_{1, 1}(\rho_k)\}\|_{l_2} \leq [1]\|Q\| \leq [1]\|q\|^2 \leq [1]\|\ff\|^2.
		\end{align}
		
		Вместе с \eqref{eq:J12} эта оценка влечет \eqref{eq:Jj} при $j=1.$
		
		\textbf{Шаг 2.2.} Используя равенство  \eqref{eq:Main1}, перепишем интеграл $J_2(\rho)$ в виде
		\begin{align*}
		2J_2(\rho) &= \int_0^{\pi} q(t) \sin \rho(\pi-t)
		\left(\frac{1}{\rho} \int_0^t M(s) \sin \rho(t-s) ds - \right. \\
		&\left. - \int_0^t (t-s)M(s)\cos \rho(t-s) ds\right) dt  = J_{2, 1}(\rho) + J_{2, 2}(\rho).
		\end{align*}
		
		В силу Леммы 2 функция
		\begin{align*}
		F(t) = \int_0^t M(s) \sin \rho(t-s) ds \in L_2(0, \pi), \quad \|F\| \leq [1]\|M_0\|.
		\end{align*}
		
		Тогда в силу неравенства Коши-Буняковского и оценки \eqref{eq:rho}
		\begin{align}\label{eq:J21}
		|J_{2, 1}(\rho_k)| \leq \frac{[1]}{\rho_k} \|q\|\|M_0\|, \quad \text{поэтому}\quad \|\{J_{2, 1}(\rho_k)\}\|_{l_2} \leq [1]\|\ff\|^2.
		\end{align}
		
		Далее, положим
		\begin{align*}
		\alpha(\rho, t) = \int_0^t (t-s) M(s) \cos \rho(t-s)ds = \int_0^t g(t, s) \cos \rho(t-s) ds,
		\end{align*}
		где $g(t, s) = (t-s)(\pi-s)^{-1}M_0(s).$ Так как $|g(t, s)|\leq |M_0(s)|,$ то в силу Леммы 4 последовательность $\{\alpha(\rho_k, t)\}_{1}^{\infty} \in l_2,$ причем
		\begin{align*}
		\|\{\alpha(\rho_k, t)\}\|_{l_2} \leq [1]\|M_0\|.
		\end{align*}
		
		Тогда из Леммы 5 следует
		\begin{align*}
		\|\{J_{2, 2}(\rho_k)\}\|_{l_2} \leq [1] \|M_0\|\|q\|\leq [1]\|\ff\|^2.
		\end{align*}
		
		Эта оценка вместе с \eqref{eq:J21} влечет \eqref{eq:Jj} при $j=2.$
		
		\textbf{Шаг 2.3.} Оценка интеграла $J_3(\rho)$ сводится к оценке  $J_2(\rho).$ Переставим в интеграле $J_3(\rho)$ порядок интегрирования с учетом $0\leq \xi \leq s \leq t \leq \pi.$ Тогда получим
		\begin{align*}
		&J_3(\rho) = \int_0^{\pi} q(\xi) \sin \rho \xi \int_{\xi}^{\pi} \sin \rho(s-\xi) \int_s^{\pi} M(t-s) \sin \rho(\pi-t) dt ds d\xi.
		\end{align*}
		
		Делая последовательно замены $t\to \pi-t,$ $s\to \pi-s,$ $\xi\to \pi-\xi,$ преобразуем этот интеграл к виду
		\begin{align*}
		&J_3(\rho) = \int_0^{\pi} q(\xi) \sin \rho \xi \int_{\xi}^{\pi} \sin \rho(s-\xi) \int_0^{\pi-s} M(\pi-t-s) \sin \rho t dt ds d\xi = \\
		&= \int_0^{\pi} q(\xi) \sin \rho \xi \int_0^{\pi-\xi} \sin \rho(\pi-s-\xi) \int_0^s M(s-t)\sin \rho t dt ds d\xi = \\
		&=\int_0^{\pi} q(\pi-\xi) \sin \rho(\pi-\xi) \int_0^{\xi} \sin\rho(\xi-s) \int_0^s M(s-t)\sin \rho t dt ds d\xi.
		\end{align*}
		
		Тем самым, $J_3(\rho)$ совпадает с $J_2(\rho),$ если в представлении для $J_2(\rho)$ функцию $q(\xi)$ заменить на $q(\pi-\xi).$ Поэтому оценка для $J_3(\rho)$ проводится так же, как для $J_2(\rho).$
		
		\textbf{Шаг 2.4.} Интеграл $J_4(\rho)$ перепишем с учетом равенства \eqref{eq:Main1}
		\begin{align}\label{eq:J4}
		\begin{aligned}
		2J_4(\rho) &= \int_0^{\pi} \sin \rho(\pi-t) \int_0^t M(t-s) \left(\frac{1}{\rho}\int_0^s M(\xi)
		\sin \rho(s-\xi) d\xi -\right. \\
		&\left.-\int_0^s (s-\xi) M(\xi) \cos \rho(s-\xi) d\xi   \right) ds dt.
		\end{aligned}
		\end{align}
		
		Обозначим
		\begin{align*}
		F(s, \rho) = \int_0^s M(\xi) \sin \rho(s-\xi) d\xi, \quad \alpha(\rho, s) = \int_0^s \frac{(s-\xi)}{\pi-\xi} M_0(\xi) \cos \rho(s-\xi)d\xi.
		\end{align*}
		
		Меняя в \eqref{eq:J4} порядок интегрирования, получаем
		\begin{align*}
		2J_4(\rho) = \int_0^{\pi} \int_s^{\pi} M(t-s) \sin \rho(\pi-t) dt \left(\frac{1}{\rho}F(s, \rho) - \alpha(\rho, s)\right) ds = J_{4, 1}(\rho) + J_{4, 2}(\rho).	
		\end{align*}
		
		В силу Леммы 2 функция
		\begin{align*}
		F_1(s, \rho) = \int_s^{\pi} M(t-s) \sin \rho(\pi-t) dt \in L_2(0, \pi), \quad \|F_1\|\leq [1]\|M_0\|.
		\end{align*}
		Функция $F(s, \rho)$ в силу Леммы 2 допускает такую же оценку $\|F\|\leq [1]\|M_0\|.$ Тогда
		\begin{align*}
		|J_{4, 1}(\rho_k)|\leq \frac{1}{|\rho_k|}\|M_0\|^2, \quad \|\{J_{4, 1}(\rho_k)\}\|_{l_2} \leq [1]\|M_0\|^2.
		\end{align*}
		
		В силу Леммы 4 для $l_2$-нормы последовательности $\{\alpha(\rho_k, s)\}_1^{\infty}$ верна оценка $\leq [1]\|M_0\|.$ Но тогда из Леммы 5 получаем оценку
		\begin{align*}
		\|\{J_{4, 2}(\rho_k)\}\|_{l_2} \leq [1]\|M_0\|^2.
		\end{align*}
		Тем самым оценка \eqref{eq:Jj} верна при $j=4.$

		\textbf{Шаг 3.} Обозначим
		\begin{align}\label{eq:R}
		R(\rho) = \rho\sum_{n=3}^{\infty} S_n(\pi, \rho).
		\end{align}
		Оценку функции $R$ проведем так же, как в \eqref{eq:Phi}. Учитывая, что $\nu = |Im \rho| \leq 1/2,$ получаем
		\begin{align}\label{Rrho}
		|R(\rho)| \leq \frac{e^{\nu \pi}}{\sqrt{3!}}t^3\left(1+\frac{t}{\sqrt{2}} + \left(\frac{t}{\sqrt{2}}\right)^2+...\right) \leq \frac{[1]t^3}{\sqrt{2}-t}, \quad t = \frac{\sqrt{\pi}\|\ff\|}{|\rho|}.
		\end{align}
		В наших предположениях $\|\ff\| <\sqrt\pi/32$,  поэтому  $\|\ff\|^3 < \|\ff\|^2.$ Тогда из  оценок \eqref{Rrho} и \eqref{eq:rho}   получаем
		\begin{align}\label{eq:Rest}
		|R(\rho_k)| \leq \frac{[1]\|\ff\|^3}{k^3} \leq \frac{1}{k^2} \beta_k, \quad \beta_k \leq [1]\|\ff\|^2 k^{-1}.
		\end{align}
		\noindent
		Cледовательно,  $l_2-$норма последовательности $\{\beta_k\}_1^{\infty}$ оценивается величиной ${\cal O}(\ff^{\,2}).$
		
		\textbf{Шаг 4.} Будем искать асимптотику нулей $\rho_k$ в виде \eqref{eq:lak}, где $\{\delta_k\}_1^{\infty}$ --- последовательность чисел, о которой, согласно Теореме 3, мы имеем априорную информацию
		\begin{align}\label{eq:apri}
		\delta_k k^{-1} = {\cal O}(\ff).
		\end{align}
		
		Перепишем уравнение
		\begin{align*}
		\rho^{-1}\sin \rho \pi + \rho^{-1}\Phi(\rho) = 0
		\end{align*}
		с учетом обозначений \eqref{eq:I}, \eqref{eq:J} и \eqref{eq:R} в виде
		\begin{align}\label{eq:equa}
		\sin \rho \pi + \rho^{-1} I(\rho) + \rho^{-2} J(\rho) + R(\rho) = 0.
		\end{align}
		Подставим
		\begin{align}
		\rho = k +r_k, \quad r_k = \frac{[q]}{2k} - \frac{a_{2k}}{4k} + \frac{m_{2k}}{4k} -\frac{n_k}{4k} +\frac{\delta_k}{k^2},
		\end{align}
		в это уравнение и воспользуемся оценками \eqref{eq:asy}, \eqref{eq:Jj}, \eqref{eq:Rest} и \eqref{eq:rho}. Из этих оценок следует
		\begin{align*}
		\rho_k^{-1} I(\rho_k) + \rho_k^{-2} J(\rho_k) + R(\rho_k) &= (-1)^{k+1} \pi(r_k-\delta_k k^{-2}) + k^{-1}\{k^{-1}\}\alpha_k +k^{-2}\beta_k = \\
		&= (-1)^{k+1} \pi r_k + k^{-2} (\beta_k+\beta'_k),
		\end{align*}
		где  $\|\{\beta_k\}\|_{l_2} = {\cal O}(\ff^{\,2}),$ $\|\{\alpha_k\}\|_{l_2} = {\cal O}(\ff).$
		Согласно обозначениям \eqref{eq:abre} тогда  $k^{-1}\{k^{-1}\}\alpha_k = k^{-2}\beta'_k$,  где
		$\beta'_k = {\cal O}(\ff^{\,2})$.

		Имеем также
		\begin{align*}
		\sin (k+r_k)\pi = \cos k\pi \sin r_k \pi = (-1)^k  r_k + {\cal O}(r_k^3).
		\end{align*}
		
		Поэтому из уравнения  \eqref{eq:equa} получаем
		\begin{align}
		(-1)^{k+1} \pi (\delta_k k^{-2}) + k^{-2} (\beta_k+\beta'_k) + {\cal O}(r_k^3) = 0.
		\end{align}
		
		Но согласно \eqref{eq:apri} имеем $r_k^3 = {\cal O}(\ff^{\,3})k^{-3} = {\cal O}(\ff^{\,2})k^{-3}.$ Поэтому из последнего уравнения получаем
		\begin{align}
		\delta_k = (-1)^{k+1}\pi^{-1}(\beta_k+\beta_k') +{\cal O}(\ff^{\,2})k^{-1},
		\end{align}
		что влечет оценку \eqref{eq:f2}.

		\textbf{Шаг 5.} Теорема доказана при условии, что наименьшее число $k_1,$ при котором выполнена оценка \eqref{eq:c0}, равно 1. Если же $k_1>1,$ то без изменений получаем
		\begin{align} \label{eq:k0'}
		\left(\sum_{k=k_1}^{\infty} |\delta_k|^{2} \right)^{1/2} = {\cal O}(\ff^{\,2}).
		\end{align}
		
		Число $k_1$  было определено при доказательстве Теоремы 3 неравенствами
		\begin{equation}\label{eq:k1}
		c_0\|\ff\| -1 \leqslant k_1 - 1 <  c_0\|\ff\|, \quad c_0 = 32/\sqrt\pi.
		\end{equation}
		Положим
		$$
		\tau_k = 2[q] -a_{2k} +m_{2k} -  n_k, \ \, \text{ тогда } \, \  \delta_k = \left(\rho_k - k - \frac{\tau_k}{4k}\right) k^2.
		$$
		Очевидно,
		$$
		|\tau_k| \leq 2\|q\| +\|q\| + 4 \|M_0\| \leq 3 \|\ff\|.
		$$
		Для оценки суммы $\sum |\delta_k|^2$  в пределах от 1 до $k_1 -1$   воспользуемся оценкой \eqref{eq:c3}. Тогда
		\begin{align*}
		&\sum |\delta_k|^2 \leq \sum 2\left(|\rho_k - k|^2 + \frac{|\tau_k|^2}{16k^2}\right) k^4 \leq
		2\left(c_2^2 \|\ff\|^4 +\frac{9}{16}\|\ff\|^2\right) \sum k^2  \leq  \\
		&\leq \left(c_2^2 \|\ff\|^4 + \|\ff\|^2\right) k_1^3   \leq \left(c_2^2 \|\ff\|^4 + \|\ff\|^2\right)8c_0^3\|\ff\|^3 \leq 8c_0^3(c_2^2 r^3 + r) \|\ff\|^4,
		\end{align*}
		если $\|\ff\|\leq r.$
		При переходах в последних неравенствах  мы воспользовались оценкой \\ $\sum k^2 \leq \int_1^{k_1}t^2\, dt < k_1^3/3$
		и оценкой  \eqref{eq:k1},  из которой следует $k_1 \leq 2c_0\|\ff\|$, если $k_1 > 1$.   Тем самым, с учетом \eqref{eq:k0'}  получаем
		$$
		\|\{\delta_k\}\|_{l_2} \leq (C +c_4 r^{3/2})\|\ff\|^2,
		$$
		если $\ff$  лежит в шаре радиуса $r$.
		Теорема доказана.

		\section{Доказательство Теоремы 2.}
		
		Пусть $\mu_k$ --- собственные значения Дирихле-Неймана. Очевидно, числа $\rho'_k = \mu_k^{1/2}$  совпадают с нулями функции $S'(\pi, \rho),$ где
		\begin{align}\label{3:eq:90}
		\begin{aligned}
		S'(x, \rho) &= \sum_{n=0}^{\infty}S'_n (x,\rho), \quad  S_0'(x,\rho) = \cos(\rho x), \\
		S'_n(x, \rho) &= \int_0^x \cos \rho (x-t) \left[ q(t) S_{n-1} (t,\rho) + \int_0^t M(t-s) S_{n-1} (s, \rho) ds \right] dt.
		\end{aligned}
		\end{align}
		
		Так же, как в Лемме 3, легко получить оценки
		\begin{align}\label{3:eq:D_n}
		|S'_n(x,\rho)| \leq \frac{\|\ff\|^n e^{\nu x} x^{n/2}}{\sqrt{n!}|\rho|^n}, \quad n\geq 0.
		\end{align}
		
		Функцию $S'(\pi, \rho)$ представим в виде
		\begin{align}
		S'(\pi, \rho) = \cos \rho \pi + \Psi(\rho).
		\end{align}
		Из оценок \eqref{3:eq:D_n}  следует (см. доказательство Теоремы 3)	
		\begin{align}
		|\Psi(\rho)| \leq \sum_{n=1}^{\infty} |S'_n(\pi, \rho)| \leq \frac{e^{\nu \pi} \sqrt{2}t}{\sqrt{2}-t}, \quad \nu = |Im \rho|, \quad t = \frac{\sqrt{\pi}\|\ff\|}{|\rho|}.
		\end{align}
		Далее, все рассуждения в доказательстве Теоремы 3 сохраняются с единственнным изменением, что
		окружности
		$$
		\gamma_k =\{\rho: |\rho - k| = r_k\} \ \, \text{надо заменить на} \ \, \gamma'_k =\{\rho: |\rho - k+1/2| = r_k\}.
		$$
		Функция $\cos \rho \pi$ на  окружностях $\gamma'_k$ оценивается снизу в точности той же величиной, что и $|\sin \rho \pi|$ в задаче Дирихле:
		\begin{align*}
		\min_{|\rho-k-1/2| = r} |\cos \rho \pi| = \min_{|\rho-k| = r} |\sin (\rho-k)\pi| \geq \pi r \left(1 - \frac{10\pi^2 r^2}{3(20-\pi^2r^2)} \right).
		\end{align*}
		Аналогично (повторяем рассуждения в Теореме 3), функция  $\Psi(\rho)$ на окружностях $\gamma'_k$ допускает оценку сверху
		\begin{align*}
		|\Psi(\rho)| < \frac{e^{\pi r}\sqrt{\pi}\|\ff\|}{k-\frac 12 - r - \sqrt{\frac{\pi}{2}} \|\ff\|}, \quad 0<r\leq \frac{1}{2}.
		\end{align*}
		Теперь, снова повторяя рассуждения в доказательстве Теоремы 3 (нужно только вместо числа $k$  писать  $k-1/2$),  придем к оценкам \eqref{eq:c0}   и \eqref{eq:c1}  с теми же постоянными
		$c_0$  и $c_1$.  В частности,
		\begin{equation}\label{eq:rhonew}
		|\rho_k' - k+1/2|< c_1\|\ff\|(k-1/2)^{-1}.
		\end{equation}
		
		Дальнейшее доказательство, как и ранее, разобьем на 5 шагов. Подробно остановимся только на наиболее важном первом шаге --- вычислении функции  $S'_1(\pi, \rho).$
		
		Из асимптотики \eqref{eq:rhonew}  получаем представления
		\begin{align}
		&\rho_k' = k - 1/2 + \{k^{-1}\}, \nonumber \\
		&\cos \rho_k' t  = \cos (k-1/2)t +\{k^{-1}\}t \sin (k-1/2)t + [1]\{k^{-2}\}, \label{eq:cos'}\\
		&\sin \rho_k' t = \sin (k-1/2)t + \{k^{-1}\} t \cos (k-1/2)t + [1]\{k^{-2}\}, \label{eq:sin'}\\
		&\frac{1}{\rho_k'}  = \frac{1}{k-1/2} + \frac{1}{k^2} \{k^{-1}\}, \quad k\geqslant k_1. \label{eq:rho1}
		\end{align}
		
		Из равенства  \eqref{eq:S1t} получаем
		\begin{align}\label{eq:5.10}
		\begin{aligned}
		S_1'(t, \rho) & = \frac{1}{2\rho}\left(\sin \rho t \int_0^t q(\xi) d\xi -   \sin \rho t \int_0^t q(\xi) \cos 2\rho \xi d\xi + \right.\\
		&+ \left. \cos \rho t \int_0^t q(\xi) \sin 2\rho \xi d\xi +  \int_0^t (t-\xi)M(\xi) \sin \rho(t-\xi) d\xi \right).
		\end{aligned}
		\end{align}
		Поясним, что после дифференцирования все  слагаемые без интегралов, а также  все интегральные  слагаемые с функцией $M$, кроме одного, сокращаются. Тогда
		$$
		S'_1(\pi, \rho) = K_1(\rho) +  K_2(\rho) +  K_3(\rho) +  K_4(\rho),
		$$
		где $K_j(\rho)$ интегральные слагаемые при  $t=\pi$ в правой части \eqref{eq:5.10}.
		
		С учетом равенств \eqref{eq:cos'}, \eqref{eq:sin'}, \eqref{eq:rho1} получаем
		\begin{align*}
		K_1 &= \frac{\sin (k-1/2)\pi}{2(k-1/2)}\pi [q] +k^{-1}\{k^{-2}\} = -\frac{(-1)^k}{2k-1}\pi [q] +k^{-1}\{k^{-2}\}, \\
		K_2 &= - \frac{\sin (k-1/2)\pi}{2(k-1/2)}   \int_0^{\pi} q(\xi) \cos 2(k-1/2) \xi d\xi +k^{-1}\{k^{-2}\} = \\
		&= \frac{(-1)^k\pi}{2(2k-1)}  a_{2k-1} +k^{-1}\{k^{-2}\}, \\
		K_3 &= \frac{\{k^{-1}\} \pi \sin (k-1/2)\pi }{2k-1} \int_0^{\pi} q(\xi) \sin (2k-1)\xi d\xi + k^{-1} \{k^{-2}\} = k^{-1}\{k^{-1}\}\alpha_k, \\
		K_4 &= \frac{1}{2k-1} \int_0^{\pi} M_0(\xi) \sin (k-1/2)(\pi-\xi) d\xi + k^{-1}\{k^{-1}\} \alpha_k =  \\
		&=\frac{\sin (k-1/2)\pi}{2k-1} \int_0^{\pi} M_0(\xi)  \cos (k-1/2) \xi d\xi +  k^{-1}\{k^{-1}\} \alpha_k  = \\
		&= \frac{(-1)^{k+1}\pi}{2(2k-1)}  m_{2k-1} + k^{-1}\{k^{-1}\} \alpha_k.
		\end{align*}	
		Складывая все 4 слагаемых, получаем
		\begin{align}\label{eq:5.15}
		S'_1(\pi, \rho_k') = (-1)^k \pi \left(-\frac{[q]}{2k-1}+\frac{a_{2k-1}}{2(2k-1)} - \frac{m_{2k-1}}{2(2k-1)}\right) + k^{-1}\{k^{-1}\}\alpha_k,
		\end{align}
		где последовательность $\{\alpha_k\}_1^{\infty}$ такова, что её $l_2$-норма подчинена оценке \eqref{eq:alphak}.
		
		На шаге 2  функцию $L(\rho) =\rho^2 S'_2(\pi,  \rho)$, как и ранее, представляем в виде четырех интегралов.
		Все эти интегралы  оцениваются величиной $[1] \|\ff\|^2$ также, как на шаге 2 в Теореме 1.
		В результате получим
		\begin{align}\label{eq:Lj}
		\|\{L(\rho_k')\}\|_{l_2} \leq [1]\|\ff\|^2.
		\end{align}
		
		На шаге 3 для функции $R(\rho) = \sum_{n\geqslant 3} S'_n(\pi, \rho)$  имеем очевидную оценку
		\begin{align}\label{eq:R'}
		|R(\rho'_k)| \leq [1] \|\ff\| k^{-3},\ \  \text{ что влечет} \ \ \|\{k^2 R(\rho'_k)\}\|_{l_2} \leq [1] \|\ff\|.
		\end{align}

		На шаге 4 уравнение $\cos \rho \pi + \Psi(\rho) = 0$ перепишем в виде
		\begin{align}\label{eq:equa'}
		\cos \rho \pi + S'_1(\pi, \rho) + \rho^{-2} L(\rho) + R(\rho) = 0.
		\end{align}
		Асимптотику нулей $\rho_k'$ ищем в виде \eqref{eq:muk}, где
		$\{\delta_k\}_1^{\infty}$ --- последовательность чисел, о которой, согласно Теореме 3',  имеется априорная информация
		\begin{align}\label{eq:apri'}
		\delta_k k^{-1} = {\cal O}(\ff).
		\end{align}
		Подставим
		\begin{align}\label{eq:rho'}
		\rho = k -\frac{1}{2} +r_k, \quad r_k = \frac{[q]}{2k-1} - \frac{a_{2k-1}}{2(2k-1)} + \frac{m_{2k-1}}{2(2k-1)} +\frac{\delta_k}{k^2},
		\end{align}
		в  уравнение \eqref{eq:equa'} и воспользуемся оценками  \eqref{eq:5.15}, \eqref{eq:Lj}, \eqref{eq:R'} и \eqref{eq:rho1}. В результате получим
		\begin{align*}
		S_1'(\pi, \rho_k') + (\rho_k')^{-2} L(\rho_k') + R(\rho_k') &= (-1)^{k+1} \pi(r_k-\delta_k k^{-2}) + k^{-1}\{k^{-1}\}\alpha_k +k^{-2}\beta_k = \\
		&= (-1)^{k+1} \pi r_k + k^{-2} \{\beta_k\},
		\end{align*}
		где
		$\|\{\beta_k\}\|_{l_2} = {\cal O}(\ff^{\,2}), \ \|\{\alpha_k\}\|_{l_2} = {\cal O}(\ff)$.
		Тогда согласно обозначениям \eqref{eq:abre} \\ $k^{-1}\{k^{-1}\}\alpha_k = k^{-2}\beta'_k$,  где
		$\beta'_k =  {\cal O}(\ff^{\,2})$.
		
		Имеем также
		\begin{align*}
		\cos (k- 1/2 +r_k)\pi = -\sin (k-1/2)\pi \sin r_k \pi = (-1)^k  r_k \pi + {\cal O}(r_k^3).
		\end{align*}
		Поэтому из уравнения  \eqref{eq:equa'} получаем
		\begin{align}
		(-1)^{k+1} \pi (\delta_k k^{-2}) + k^{-2} \beta_k + {\cal O}(r_k^3) = 0.
		\end{align}	
		Но согласно \eqref{eq:rho'}  и \eqref{eq:apri'} имеем $r_k = k^{-1} {\cal O}(\ff)$,  поэтому
		$r_k^3 = {\cal O}(\ff^{\,3})k^{-3} = {\cal O}(\ff^{\,2})k^{-3}.$  Тогда из последнего уравнения получаем
		\begin{align}
		\delta_k = (-1)^{k+1}\pi^{-1}\beta_k +{\cal O}(\ff^{\,2})k^{-1}.
		\end{align}
		что влечет оценку \eqref{eq:f2} для  асимптотики собственных значений задачи  Дирихле-Неймана.
		
		Доказательство на шаге 5 в случае $k_1 > 1$ проводится без изменений. Этим завершается доказательство
		Теоремы 2.

	
\end{document}